\documentclass[a4paper,reqno,12pt]{amsart}

\usepackage{amsmath}
\usepackage{amscd}
\usepackage{amssymb}
\usepackage{mathrsfs}
\usepackage{stmaryrd}

\newtheorem{theorem}{Theorem}[section]

\newtheorem{lemma}{Lemma}[section]
\newtheorem{corollary}{Corollary}[section]

\newtheorem{remark}{Remark}[section]

\newcommand{\Ric}{\mathrm{Ric}}

\newcommand{\del}{\nabla}
\newcommand{\laplace}{\triangle}

\newcommand{\tq}{\widetilde{q}}
\newcommand{\tp}{\widetilde{p}}
\def\Ric{{\rm Ric}\,}

\def\be{\begin{equation}}
\def\ee{\end{equation}}

\def\endproof{{\hfill $\square$}\medskip}

\AtEndDocument{\bigskip{\footnotesize
		\textsc{School of Mathematics (Zhuhai), Sun Yat-Sen University (Zhuhai), Zhuhai, China}\par
		\textit{E-mail address: }\texttt{maiweixiong@gmail.com}\par
		\medskip
		\textsc{Shanghai Center for Mathematical Sciences, Fudan University, Shanghai, China}\par
		\textit{E-mail address: }\texttt{oujianyu@fudan.edu.cn}
}}

\begin{document}

\title{Uncertainty Principle and its rigidity on complete gradient shrinking Ricci solitons}

\author{Weixiong Mai \& Jianyu Ou}

\maketitle

\begin{abstract}
We prove rigidity theorems for shrinking gradient Ricci solitons supporting the Heisenberg-Pauli-Weyl uncertainty principle with the sharp constant in $\mathbb{R}^n$. In addtion, we partially give analogous rigidity results of the Caffarelli-Kohn-Nirenberg inequalities on shrinking Ricci solitons.
\end{abstract}

\section{Introduction}
In quantum mechanics, the Heisenberg uncertainty principle addresses that the position and the momentum of a particle cannot be both determined precisely in any quantum state, which in mathematics is stated as:\\
\noindent{\it
	For a function $u\in C_0^\infty(\mathbb R^n),$ there holds
%	\begin{align}\label{HUP}
%	\int_{\mathbf R^n}|x|^2|u(x)|^2 dx\int_{\mathbf R}|\xi|^2 |\hat u(\xi)|^2d\xi \geq \frac{n^2}{4} \left(\int_{\mathbf R^n}|u(x)|^2dx\right)^2,
%	\end{align}
%	where $\hat u(\xi)$ is the Fourier transform of $u(x).$ Or equavilently,
		\begin{align}\label{HUP1}
	\int_{\mathbb R^n}|x|^2|u(x)|^2 dx\int_{\mathbb R^n}|\del u(x)|^2dx \geq \frac{n^2}{4} \left(\int_{\mathbb R^n}|u(x)|^2dx\right)^2.
	\end{align}
}\\
Since a rigorous mathematical formulation of the Heisenberg uncertainty principle was first given by Pauli and Weyl (see e.g. \cite{Weyl}), it is called appropriately the Heisenberg-Pauli-Weyl (HPW) uncertainty principle. In (\ref{HUP1}) HPW uncertainty principle is of the form in terms of gradients, which could be generalized to general domains and manifolds, while in literatures HPW uncertainty principle is originally given by using the Fourier transform. 
HPW uncertainty principle has been developed a lot due to its significance in both physics and mathematics. In the Euclidean space, the developments of HPW uncertainty principle include, for instance, $L^p$-version of HPW uncertainty principle, the study of HPW uncertainty principle in signal analysis, and an abstract version of uncertainty principle for two self adjoint operators on Hilbert spaces from the aspect of operator theory. See e.g. \cite{Fefferman,Folland-Sitaram,Cohen,Dangetal} and the references therein for details.

 In recent, HPW uncertainty principle on compact or non-compact manifolds has attracted researchers' attention. For manifold cases, the study of HPW uncertainty principle is based on (\ref{HUP1}) since the Fourier transform is not well-defined on manifolds in general. In Erb's thesis (\cite{Erb}), a systematical study of HPW uncertainty principle on Riemaniann manifolds (especially on compact manifolds) is given, which is essentially based on Breitenberger's generalization of uncertainty principle (\cite{Breitenberger}). In \cite{Kombe-Ozaydin,Kombe-Ozaydin1}, Kombe and \"Ozaydin study the sharp HPW uncertainty principle on complete manifolds. It is noted that HPW uncertainty principle can be viewed as an endpoint of the well-known Caffarelli-Kohn-Nirenberg (CKN) inequalities (\cite{Caffarelli-Kohn-Nirenberg}, and see also a brief introduction to CKN inequalities below). A rigidity result of a subclass of CKN inequalities on manifolds with non-negative Ricci curvature is proved in \cite{Xia} by Xia. 
  In \cite{Kri} Krist\' aly proves a rigidity theorem for HPW uncertainty principle on complete manifolds with non-negative Ricci curvature, and a weaker rigidity theorem for a subclass of CKN inequalities on Cartan-Hadamard manifolds, which generalize Xia's result. Here we are particularly interesting in the rigidity theorem for HPW uncertainty principle, which is stated as follows. Denote by $C_0^\infty(M)$ the function space consisting of all real-valued smooth functions with compact support on $M$.\\ 
\noindent \textbf{Theorem A}(\cite[Theorem 1.2]{Kri}) {\it Let $(M^n,g)$ be a complete, $n$-dimensional Riemannian manifold with non-negative Ricci curvature. The following statements are equivalent:\\
	(a) HPW uncertainty principle 
	\begin{align*}
	\left(\int_M d_{x_0}^2 u^2 dV\right) \left(\int_M |\del u|^2dV\right)\geq \frac{n^2}{4}\left(\int_M u^2dV\right)^2, \quad u\in C^\infty_0(M),
	\end{align*}
	holds for some $x_0\in M$;\\
	(b) HPW uncertainty principle 
	\begin{align*}
	\left(\int_M d_{x_0}^2 u^2 dV_g\right) \left(\int_M |\del u|^2dV\right)\geq \frac{n^2}{4}\left(\int_M u^2dV\right)^2, \quad u\in C^\infty_0(M),
	\end{align*}
	holds for every $x_0\in M$;\\
	(c) $(M^n,g)$ is isometric to $\mathbb R^n.$
}\\
Here $d_{x_0}$ the distance function from a point $x_0\in M$ and $dV$ is the canonical volume element.
In \cite{Nguyen} Nguyen  generalizes Krist\' aly's and Xia's results to CKN inequalities with a larger  class of parameters. To state Nguyen's results, we recall a subclass of CKN inequalities given in \cite{Nguyen}:
let $n\geq 2,p>1,r>0$ and $\alpha,\beta,\gamma\in \mathbb R$ such that 
\begin{align}\label{cond1}
\begin{split}
\frac{1}{r}-\frac{\gamma}{n}>0, \frac{1}{p}-\frac{\alpha}{n}>0,1-\frac{\beta}{n}>0
\end{split}
\end{align}
and 
\begin{align}\label{cond2}
\begin{split}
\gamma = \frac{1+\alpha}{r}+\frac{p-1}{pr}\beta.
\end{split}
\end{align}
It is well-known (\cite{Caffarelli-Kohn-Nirenberg}) that for any $u\in C_0^\infty(\mathbb R^n)$ there holds
\begin{align}\label{CKN-RN}
\left(\int_{\mathbb R^n}\frac{|\del u|^p}{|x|^{\alpha p}}dx\right)^\frac{1}{p} \left(\int_{\mathbb R^n} \frac{|u|^\frac{p(r-1)}{p-1}}{|x|^\beta}dx\right)^\frac{p-1}{p}\geq \frac{n-\gamma r}{r}\int_{\mathbb R^n}\frac{|u|^r}{|x|^{\gamma r}}dx.
\end{align}
Moreover, the sharpness of the constant $\frac{n-\gamma r}{r}$ and the extremal function with different ranges of parameters are specifically given in \cite{Nguyen}. As mentioned previously, the CKN inequalities (\ref{CKN-RN}) give HPW inequality when $p=r=2,\alpha=0,\beta=-2,\gamma =0.$

\noindent \textbf{Theorem B}(\cite[Theorem 1.6]{Nguyen}) {\it Given $n\geq 2,p>1,$ and $\alpha,\beta,\gamma\in \mathbb R$ such that (\ref{cond1}) and (\ref{cond2}) hold true with $r=p.$ Suppose, in addition, that $1+\alpha-\frac{\beta}{p}>0.$ Let $(M^n,g)$ be a complete, $n$-dimensional Riemannian manifold with non-negative Ricci curvature. The following statements are equivalent:\\
	(a) For any $u\in C^\infty_0(M),$ CKN inequalities 
	\begin{align*}
	\left(\int_M d_{x_0}^{-\beta} u^{\frac{p(r-1)}{p-1}} dV\right)^\frac{p-1}{p} \left(\int_M d_{x_0}^{-\alpha p}|\del u|^pdV\right)^\frac{1}{p}\geq \frac{n-\gamma r}{r}\int_M d_{x_0}^{-\gamma r}u^rdV
	\end{align*}
	holds for some $x_0\in M$;\\
	(b) For any $u\in C^\infty_0(M),$ CKN inequalities 
	\begin{align*}
\left(\int_M d_{x_0}^{-\beta} u^{\frac{p(r-1)}{p-1}} dV\right)^\frac{p-1}{p} \left(\int_M d_{x_0}^{-\alpha p}|\del u|^pdV\right)^\frac{1}{p}\geq \frac{n-\gamma r}{r}\int_M d_{x_0}^{-\gamma r}u^rdV
\end{align*}
	holds for every $x_0\in M$;\\
	(c) $(M^n,g)$ is isometric to $\mathbb R^n.$
}\\
Obviously, Theorem B contains Theorem A.
We refer to \cite{Nguyen} for the rigidity results of CKN inequalities with general parameters.
As noted, Theorem A and Theorem B can be, in fact, included into the best constant program initiated by Aubin \cite{Aubin}, and studied by Ledoux \cite{Ledoux}, Bakry, Concordet and Ledoux \cite{Bakry-Concordet-Ledoux}, Cheeger and Colding \cite{Cheeger-Colding}, Druet, Hebey and Vaugon \cite{Druet-Hebey-Vaugon}, do Carmo and Xia \cite{Carmo-Xia}, Minerbe \cite{Minerbe}, Li and Wang \cite{Li-Wang}, Xia \cite{Xia1,Xia2,Xia}, Krist\'aly \cite{Kri,Kri1}, Krist\'aly and Ohta \cite{Kri-Ohta}, Nguyen \cite{Nguyen}, etc.

It is natural and of significance to ask if there exist some rigidity properties of HPW uncertainty principle (or more general CKN inequalities) on the complete manifold with non-negative Bakry-Emery Ricci curvature, i.e., $\Ric_f\geq 0,$ where $\Ric_f=\Ric+\del\del f$ for some smooth function $f$ on $M.$ As a special case, gradient shrinking Ricci solitons are with non-negative Bakry-Emery Ricci curvature. 

A complete Riemannian manifold $(M^n, g_{ij}, f)$ is called a gradient Ricci soliton if there exists a smooth function $f$ on $M^n$ such that the Ricci tensor $R_{ij}$ of the metric $g_{ij}$ satisfies the equation
$$R_{ij}+\del_i \del_j f=\rho g_{ij}$$
for some constant $\rho$. Here, $\del_i \del_j f$ denotes the Hessian of $f$. The Ricci soliton is called shrinking, steady, or expanding if $\rho>0,\ \rho=0$ or $\rho<0$ respectively. The function $f$ is called a potential function of the gradient Ricci soliton. It is clear that if $f$ is a constant function, the gradient Ricci soliton is simply an Einstein manifold. Thus, Ricci solitons are natural extensions of Einstein metrics. On the other hand, gradient Ricci solitons are also self-similar solutions to Hamilton's Ricci flow and play an important role in the study of formation of singularities in the Ricci flow, see \cite{Cao} for a nice survey on the subject. 

In this paper we consider analogous results on shrinking Ricci solitons, which often arise as possible Type I singularity models in the Ricci flow. In the following, we will normalize the soliton constant $\rho={1\over2}$ so that the shrinking gradient soliton equation is given by 
\begin{eqnarray}{\label{shrinker}}
R_{ij}+\del_i\del_j f={1\over2}g_{ij}.
\end{eqnarray}
We can write (\ref{shrinker}) as $Ric_f={1\over2}g$, so that shrinking Ricci solitons are non-negative Bakry-Emery Ricci curvature.
 The Bishop-Gromov volume comparison theorem (see e.g. \cite{Carmo}) plays an essential role in the proofs of rigidity theorems (like Theorem A and Theorem B) on manifolds with non-negative Ricci curvature (see e.g. \cite{Kri,Nguyen,Xia}). However, since the techniques in literatures are invalid in our case (For example, FIK shrinking Ricci soliton \cite{FIK} has a conelike end but it is not $\mathbb{R}^n$), we have to find new methods. In this paper we prove several rigidity theorems of HPW uncertainty principle and CKN inequalities on shrinking Ricci solitons by elementary methods.

\begin{theorem}{\label{UP}}
	Let $(M^n, g, f)$ be a complete non-compact gradient shrinking Ricci soliton. $R$ is the scalar curvature of $M$. Denote $\rho(x)=2\sqrt{f(x)}$, then 

	(a) The Heisenberg-Pauli-Weyl uncertainty principle
	\begin{align*}
		\int_M |\del u|^2 dV \int_M \rho^2u^2 dV \geq & {n^2\over4}\left(\int_M u^2dV\right)^2-{n}\int_M u^2 dV\int_M u^2R dV\\
		&+4\int_M u^2R\int_M|\del u|^2 dV+\left(\int_M u^2R dV\right)^2 ,
	\end{align*}
	holds for all $u\in C^{\infty}_{0}(M)$, the equality holds if and only if $u=e^{-f}$.

	(b) In particular, the inequality holds 
	$$\int_M |\del u|^2 dV \int_M \rho^2u^2 dV \geq {n^2\over4}\left(\int_M u^2dV\right)^2,$$
	if and only if $M^n$ is isometric to $\mathbb{R}^n$.
\end{theorem}

When $(M^n, g, f)$ is a complete non-compact shrinker with constant scalar curvature, by the result in {\cite{FG}}, we know $R\equiv\frac{K}{2},$ where $K\leq n-1$ is a non-negative integer. By Theorem \ref{UP} we have
\begin{align*}
\int_M |\del u|^2 dV \int_M \rho^2u^2 dV \geq & \left(\frac{n-K}{2}\right)^2 \left(\int_M u^2 dV\right)^2\\
& +{2K}\int_M u^2 dV\int_M|\del u|^2 dV.
\end{align*}
This case is interesting since it contains cylinders. To the authors' best knowledge, this HPW uncertainty principle is not included in literatures (in some sense in \cite{Erb} Erb also gives HPW uncertainty principle on cylinders, which is, however, essentially different from our case). From the above inequality, we can see exactly the relation between the scalar curvature and the constant of HPW uncertainty principle on cylinders. For this case we have the following result.
\begin{theorem}{\label{R=K/2}}
	Let $(M^n,g,f)$ be a complete non-compact gradient shrinking Ricci soliton. Denote $\rho(x)=2\sqrt{f(x)}$. For any $u\in C^{\infty}_{0}(M)$, there exists a constant $0\leq K<n$ such that the Heisenberg-Pauli-Weyl uncertainty principle inequality
	\begin{align}{\label{constant}}
		\int_M|\del u|^2dV \int_M \rho^2u^2 dV &\geq \left({n-K\over 2}\right)^2\left(\int_M u^2dV\right)^2\nonumber\\
		&+{2K}\int_M u^2 dV \int_M |\del u|^2 dV,
	\end{align}
	holds then the scalar curvature holds ${\int_M Re^{-2f}dV\over{\int_M e^{-2f}dV}}={K\over2}$.
\end{theorem}
\begin{remark}
	Let $K=0$ in the Theorem {\ref{R=K/2}}, we can also get the part (b) of Theorem {\ref{UP}}.
\end{remark}
\begin{remark}
	In \cite{CaoZhou}, we know that for any shrinker $(M^n, g, f)$ the average scalar curvature $R$ over $D(r):=\{x\in M| f(x)\leq r\}$ is bounded by ${n\over2}$. That means ${\int_{D(r)}RdV\over{\int_{D(r)}dV}}< {n\over2}$. Theorem {\ref{R=K/2}} gives us the average scalar curvature $R$ over the whole space with the weight $e^{-2f}$ when ({\ref{constant}}) holds, then ${\int_M Re^{-2f}dV\over{\int_Me^{-2f}dV}}= {K\over2}< {n\over2}$. In addition, if $R\leq \frac{K}{2}<\frac{n}{2}$ or $ \frac{K}{2}\leq R<\frac{n}{2},$ one can easily conclude that $R\equiv \frac{K}{2}$ from the proof of Theorem \ref{R=K/2}.
\end{remark}

In the following we present CKN inequalities on shrinking Ricci solitons 
\begin{theorem}\label{pigola}
	Let $(M^n,g,f)$ be a complete gradient shrinking Ricci soliton.
	Given $n\geq 2,$ $p>1,r>1$ and $\alpha,\beta,\gamma\in \mathbb R$ such that (\ref{cond1}) and (\ref{cond2}) hold true.  There holds 
	\begin{align}\label{CKN-srs}
	\begin{split}
	&\left(\int_M f^{-\frac{\alpha p}{2}}|\del u|^pdV\right)^\frac{1}{p} \left(\int_M f^{-\frac{\beta}{2}}\left(1-\frac{R}{f}\right)^\frac{p}{2(p-1)}u^\frac{p(r-1)}{p-1}dV\right)^{\frac{p-1}{p}}\\
	&\geq\left|\frac{n-\gamma r}{2r}\int_M u^r f^{-\frac{\gamma r}{2}} dV-\frac{1}{r}\int_M u^rf^{-\frac{\gamma r}{2}}RdV+\frac{\gamma r}{2}\int_M f^{-\frac{\gamma r}{2}-1}u^r RdV\right|
	\end{split}
	\end{align}
	for any $u\in C_0^\infty(M).$\\
	In particular, when the parameters satisfy
	$\alpha=0,r=\widetilde p,\gamma r={\widetilde q},p=2,\beta=2\widetilde q-2,$ where $0<\widetilde q<2<\widetilde p,\frac{2(\widetilde p-\widetilde q)}{\widetilde p-2}>n>2,$ the inequality (\ref{CKN-srs}) becomes
    \begin{align}\label{GHUP}
    \begin{split}
     &\int_M |\nabla u|^2 dV \int_M f^{-\tq+1}u^{2\tp-2}dV \\
     &\geq \left(\frac{n-\tq}{2\tp}\int_M f^{-\frac{\tq}{2}}u^{\tp} dV - \frac{1}{\tp}\int_M f^{-\frac{\tq}{2}}u^{\tp} R dV+\frac{\tq}{2\tp}\int_M f^{-\frac{\tq}{2}-1}u^{\tp} R dV
     \right)^2\\
     & + \int_M f^{-\tq}u^{2\tp-2}RdV \int_M |\nabla u|^2 dV,
     \end{split}
    \end{align}
    and moreover, the equality holds $u_\lambda=(\lambda+f^{-\frac{\tq}{2}+1})^\frac{1}{2-\tp},$ $\lambda>0.$

\end{theorem}

As an application, we can get the following corollary which is contained in \cite[Theorem 3]{PRS}.
\begin{corollary}
There doesn't exist a complete non-compact gradient shrinking Ricci soliton with $R\geq {n\over2}$.
\end{corollary}

\begin{remark}
	When $R=0,$ by Theorem \ref{pigola} we just obtain CKN inequalities on $\mathbb R^n.$ In the Euclidean space, the equalities of CKN inequalities are respectively given by different classes of functions according to different ranges of the parameters (see \cite{Nguyen}). However, this usually is not true in our case since we replace the distance function by the double square root of the potential function $\rho=2\sqrt f$.
\end{remark}
 In the following we partially give a converse result.
\begin{theorem}\label{CKN-thm}
Let $(M^n,g,f)$ be a complete gradient shriking Ricci soliton.
	Given $n\geq 2,p>1$ and $\alpha,\beta,\gamma\in \mathbb R$ such that (\ref{cond1}) and (\ref{cond2}) hold true with $r=p.$ Suppose, in addition, that $\gamma\leq 0,\ 0<1+\alpha-\frac{\beta}{p}\leq 2.$  If there holds
	\begin{align}\label{CKN1}
	\left(\int_{M} f^{-\frac{\alpha p}{2}}|\del u|^pdV\right)^\frac{1}{p} \left(\int_{M} f^{-\frac{\beta}{2}}u^p dV\right)^\frac{p-1}{p}\geq \frac{n-\gamma p}{2p}\int_{M} f^{-\frac{\gamma p}{2}}u^pdV
	\end{align}
	for any $u\in C_0^\infty(M)$, then  $M$ is isometric to $\mathbb R^n.$
\end{theorem}

\begin{remark}
	For CKN inequalities on shrinking Ricci soliton, the integrability in Theorem \ref{pigola} and Theorem \ref{CKN-thm} would not be a risk. In fact, for  $R\equiv 0,$ the integrability is guaranteed by the Euclidean case, while for $R\not\equiv 0,$ we have that $f$ is strictly positive due to the fact that $fR>\epsilon>0$ (see \cite{Chow-Lu-Yang}), where $\epsilon$ is a constant.
\end{remark}

We will prove Theorem {\ref{UP}} and Theorem {\ref{R=K/2}} in Section 3. Theorem {\ref{pigola}} and Theorem {\ref{CKN-thm}} will be proved in Section 4.

{\bf Acknowledgement} The authors thanks Prof. Huai-Dong Cao for motivating the authors to consider Theorem {\ref{R=K/2}} and Dr. Yashan Zhang and Mr. Ronggang Li for helpful discussions.

\section{Preliminaries}
In this section, we briefly recall some basic facts about gradient shrinking  Ricci solitons. Throughout the rest of the paper, we denote by
$$Rm=\{R_{ijkl}\},\ \Ric=\{R_{ij}\},\ R$$
the Riemann curvature, the Ricci tensor, and the scalar curvature of the metric $g_{ij}$, respectively.

\begin{lemma}{(Hamilton \cite{Ham1})}
Let $(M^n,g_{ij},f)$ be a complete gradient shrinking Ricci soliton satisfying the equation (\ref{shrinker}). Then, we have
\begin{eqnarray}
\del_i R=2R_{ij}\del_j f
\end{eqnarray}
and 
\begin{eqnarray}
R+|\del f|^2-f=C_0
\end{eqnarray}
for some constant $C_0$. 
\end{lemma}
Note that if we normalize $f$ by adding  the constant $C_0$, then we have
\begin{eqnarray}{\label{prop1}}
R+|\del f|^2=f
\end{eqnarray}

Now we give the following lemma which is very important in our proof.
\begin{lemma}{(B-L Chen \cite{Chen}, P-R-S \cite{PRS})}
Let $(M^n, g_{ij}, f)$ be a complete gradient shrinking Ricci soliton. Then, it has nonnegative scalar curvature. Moreover, the scalar curvature $R$ is positive, unless $(M^n, g_{ij}, f)$ is the Gaussian soliton $(\mathbb{R}^n, \delta_{ij}, {{|x|^2}\over4})$. That means if $R=0$ on a point in $M$, then $M$ is isometric to $\mathbb R^n$. 
\end{lemma}

\begin{lemma}{(Cao-Zhou \cite{CaoZhou})}
Let $(M^n,g_{ij},f)$ be a complete noncompact gradient shrinking Ricci soliton with the normalizations  (\ref{shrinker}) and (\ref{prop1}). Then,
 
(1) the potential function $f$ satisfies the estimates
\begin{eqnarray}{\label{distance}}
{1\over4}(r(x)-c_1)^2\leq f(x)\leq {1\over4}(r(x)+c_2)^2,
\end{eqnarray}
where $r(x)=d(x_0,x)$ is the distance function from a fixed point $x_0\in M$, $c_1$ and $c_2$ are positive constants depending only on the geometry of $g_{ij}$ on the unit ball $B(x_0,1)$.

(2) there exists some constant $A>0$ such that
$$Vol(B(x_0, s))\leq As^n$$
for $s>0$ sufficiently large.
\end{lemma}

\section{Uncertainty Principle on shrinking Ricci solitons}

%\begin{theorem}{\label{UP}}
%Let $(M^n, g, f)$ be a complete non-compact gradient shrinking Ricci soliton. If there holds the Heisenberg-Pauli-Weyl uncertainty principle
%\begin{eqnarray*}
%\int_M |\del u|^2 dV \int_M f|u|^2 dV &\geq& {n^2\over16}(\int_M |u|^2dV)^2-{n\over4}\int_M u^2 dV\int_M u^2R dV\\
%&&+\int_M u^2R\int_M|\del u|^2 dV+{1\over4}(\int_M u^2R)^2 dV,
%\end{eqnarray*}
%for all $u\in C^{\infty}_{0}(M)$, the equality holds if and only if $u=e^{-f}$.
%In particular, the inequality holds 
%$$\int_M |\del u|^2 dV \int_M f|u|^2 dV \geq {n^2\over16}(\int_M |u|^2dV)^2,$$
%if and only if $M^n$ is isometric to $\mathbb{R}^n$.
%
%\end{theorem}
We now prove Theorem {\ref{UP}}.

\noindent{\bf{Proof of Theorem \ref{UP}:}}

First we consider 
\begin{eqnarray*}
&&\left(\int_M 2u\sqrt{f}\langle\del \sqrt{f}, \del u\rangle dV\right)^2\\
&=&\left(\int_M u^2\sqrt{f}\laplace{\sqrt{f}}dV+\int_M u^2|\del \sqrt{f}|^2dV\right)^2\\
&=&\left(\int_M u^2\left({{n\over 4}-|\del {\sqrt{f}}|^2-{R\over2}}dV\right)+\int_M u^2|\del \sqrt{f}|^2dV \right)^2\\
&=&\left({n\over4}\int_M u^2dV-{1\over2}\int_M u^2RdV\right)^2.
\end{eqnarray*}

Using Cauchy-Schwarz's inequality twice, we have

\begin{align}\label{CS-UP}
\left(\int_M 2u\sqrt{f}\langle\del \sqrt{f}, \del u\rangle dV \right)^2\leq 4\int_M u^2 f |\del{\sqrt{f}}|^2\int |\del u|^2dV.
\end{align}

That means

\begin{eqnarray}{\label{up1}}
\int_M u^2 |\del f|^2dV\int_M |\del u|^2dV\geq \left({n\over4}\int_M u^2dV-{1\over 2}\int_M u^2 RdV\right)^2.
\end{eqnarray}

Using $|\del f|^2=f-R$, then we get 

\begin{eqnarray*}
\int_M |\del u|^2 dV \int_M f|u|^2 dV &\geq& {n^2\over16}\left(\int_M |u|^2dV\right)^2-{n\over4}\int_M u^2 dV\int_M u^2R dV\\
&&+\int_M u^2R\int_M|\del u|^2 dV+{1\over4}\left(\int_M u^2R dV\right)^2.
\end{eqnarray*}
Let $\rho=2\sqrt{f}$, we finish the proof.

Next we consider the equality in ($\ref{up1}$). Let $u=e^{-f}$,  the left-hand side in ({\ref{up1}}) is

$$\int_M u^2 |\del f|^2dV\int_M |\del u|^2dV=\left(\int_M e^{-2f}|\del f|^2dV\right)^2.$$

The right-hand side is 

\begin{eqnarray*}
&&\left({n\over4}\int_M u^2dV-{1\over 2}\int_M u^2 RdV\right)^2\\
&=&{1\over4}\left(\int_M e^{-2f} \laplace f dV\right)^2\\
&=&{1\over 4} \left(-2\int_M e^{-2f} |\del f|^2dV\right)^2\\
&=& LHS.
\end{eqnarray*}

On the other hand, we note that the equality in (\ref{up1}) holds if and only if the equality in (\ref{CS-UP}) holds. In (\ref{CS-UP}) the equality holds when $\del u=-2u\sqrt f\del\sqrt f,$ which implies that $\del \log u= -\del f.$ Hence, $u=e^{-f}.$

Finally we prove part (b). If for all $u\in C^{\infty}_0(M)$, there holds the inequality
\begin{eqnarray}
\int_M |\del u|^2 dV \int_M \rho^2u^2 dV \geq {n^2\over4}\left(\int_M u^2dV\right)^2.
\end{eqnarray}

That is 
\begin{eqnarray}{\label{upR}}
\int_M |\del u|^2 dV \int_M fu^2 dV \geq {n^2\over16}\left(\int_M u^2dV\right)^2
\end{eqnarray}
holds for $u_\lambda=e^{-\lambda f}\in \overline{C^{\infty}_0(M)}$ (see e.g. \cite[Page 49, Theorem 3.1]{Hebey}), where $0<\lambda\leq {1\over2}$ is a constant.

That means

\begin{eqnarray*}
4\lambda^2\int_M e^{-2\lambda f} |\del f|^2dV \int_M fe^{-2\lambda f}dV\geq {n^2\over4} \left(\int_M e^{-2\lambda f}dV\right)^2. 
\end{eqnarray*}

Using the shrinking soliton equation $|\del f|^2=f-R$, then

\begin{eqnarray}{\label{u}}
4\lambda^2\left(\int_M f e^{-2\lambda f}dV\right)^2&\geq& {n^2\over 4} \left(\int_M e^{-2\lambda f}dV\right)^2{\nonumber}\\
&&+4\lambda^2\int_M R e^{-2\lambda f}dV \int_M f e^{-2\lambda f}dV.
\end{eqnarray}

Let $F(\lambda)=\int_M e^{-2\lambda f}dV$, then $F'(\lambda)=-2\int_M f e^{-2\lambda f}dV<0$. So (\ref{u}) can be written as 

\begin{eqnarray*}
\lambda^2 (F'(\lambda))^2\geq {n^2\over4} (F(\lambda))^2 +4\lambda^2\int_M R e^{-2\lambda f}dV\int_M f e^{-2\lambda f}dV.
\end{eqnarray*}

Now we use contridiction method, if $M$ is not isometric to $\mathbb{R}^n$, $R>0$. Then 
$$\left({F'(\lambda)\over F(\lambda)}\right)^2>{n^2\over{4\lambda^2}}.$$

Because $F'(\lambda)<0$, we get
\begin{eqnarray}{\label{F'1}}
{F'(\lambda)\over F(\lambda)}<-{n\over 2\lambda}.
\end{eqnarray}

On the other hand, 
\begin{eqnarray*}
F'(\lambda)&=&-2\int_M f e^{-2\lambda f}dV\\
&=&-2\int_M |\del f|^2 e^{-2\lambda f}dV-2\int_M R e^{-2\lambda f}dV\\
&=&{1\over\lambda}\int_M \langle \del f, \del e^{-2\lambda f}\rangle dV -2\int_M R e^{-2\lambda f}dV\\
&=&-{1\over \lambda}\int_M {\laplace f}e^{-2\lambda f}dV-2\int_M R e^{-2\lambda f}dV\\
&=&-{1\over \lambda}\int_M \left({n\over 2}-R\right)e^{-2\lambda f}dV-2\int_M R e^{-2\lambda f}dV\\
&=&-{n\over 2\lambda} \int_M e^{-2\lambda f}dV+\left({1\over \lambda}-2\right)\int_M R e^{-2\lambda f}dV.
\end{eqnarray*}

Since $\lambda \leq {1\over 2}$, that implies
$${F'(\lambda)\over F(\lambda)}\geq -{n\over 2\lambda},$$
which is a contradiction to (\ref{F'1}).

So $R=0$ which implies $M$ is isometric to $\mathbb{R}^n$.
\endproof

%\begin{theorem}{\label{R=K/2}}
%Let $(M^n,g,f)$ be a complete non-compact gradient shrinking Ricci soliton. For any $u\in C^{\infty}_{0}(M)$, there exists a constant $K$ such that the HPW uncertainty principle type inequality holds:
%\begin{eqnarray*}%{\label{constant}}
%\int_M|\del u|^2dV \int_M fu^2 \geq ({n-K\over 4})^2(\int_M u^2dV)^2+{K\over2}\int_M u^2 dV \int_M |\del u|^2 dV,
%\end{eqnarray*}
%if and only if the scalar curvature $R\equiv {K\over2}$ is a constant. 
%\end{theorem}

Then we give the proof of Theorem \ref{R=K/2}.

\noindent{\bf{Proof of Theorem \ref{R=K/2}:}}

We also consider the function $u=e^{-f}$. Put $u=e^{-f}$ into the inequality, (\ref{constant}), 
we get
\begin{eqnarray}{\label{constant1}}
&&\int_M e^{-2f} |\del f|^2dV\int_M e^{-2f}fdV {\nonumber}
\\&\geq& \left({n-K\over 4}\right)^2\left(\int_M e^{-2f}dV\right)^2{\nonumber}\\
&&+{K\over2}\int_M e^{-2f}dV\int_M e^{-2f} |\del f|^2dV.
\end{eqnarray}

By Stokes' theorem and soliton equation $f=|\del f|^2+R$,
\begin{eqnarray}{\label{constant2}}
&&\int_M e^{-2f}|\del f|^2dV\nonumber\\
&=&-{1\over2}\int_M \langle \del f, \del e^{-2f}\rangle dV\nonumber\\
&=&{1\over2}\int_M \laplace f e^{-2f}dV\nonumber\\
&=&{1\over2}\int_M \left({n\over2}-R\right)e^{-2f}dV.
\end{eqnarray}

Using the soliton equation and (\ref{constant2}) we also have 
\begin{eqnarray}{\label{constant3}}
\int_M e^{-2f}fdV&=&\int_M e^{-2f}|\del f|^2dV+\int_M e^{-2f}RdV\nonumber\\
&=&{1\over2}\int_M \left({n\over2}-R\right)e^{-2f}dV+\int_M e^{-2f}RdV.
\end{eqnarray}

Put (\ref{constant2}) and (\ref{constant3}) into (\ref{constant1}), we obtain
\begin{eqnarray*}
\left(\int_M \left({n\over4}-{R\over2}\right) e^{-2f}dV\right)^2 &\geq& \left(\int_M \left({n-K\over4}\right)e^{-2f}dV\right)^2\\
&&+\int_M \left({K\over4}-{R\over2}\right)e^{-2f}dV\int_M \left({n\over2}-R\right)e^{-2f}dV.
\end{eqnarray*}

Consequently, we get 
\begin{eqnarray*}
&&\int_M \left({n\over2}-{R\over2}-{K\over4}\right)e^{-2f}dV\int_M \left({K\over4}-{R\over2}\right)e^{-2f}dV\\
&\geq& \int_M \left({K\over4}-{R\over2}\right)e^{-2f}dV\int_M \left({n\over2}-R\right)e^{-2f}dV.
\end{eqnarray*}

That means
$$-\left(\int_M \left({K\over4}-{R\over2}\right)e^{-2f}dV\right)^2\geq 0,$$
which implies  ${\int_M Re^{-2f}dV\over{\int_Me^{-2f}dV}}= {K\over2}$.
\endproof

\begin{remark}
	By the same method, we can get the constant $(n/4)^2$ is sharp. That means there doesn't exist a constant $K>0$ such that $$\int_M|\del u|^2dV \int_M \rho^2u^2 \geq \left({n+K\over 2}\right)^2 \left(\int_M u^2dV\right)^2-2K\int_M u^2dV\int_M |\del u|^2 dV$$
	holds for any $u\in C^\infty_0(M).$
\end{remark}

\section{CKN inequalities on shrinking Ricci solitons}
\noindent\textbf{Proof of Theorem \ref{pigola}:} For $u\in C_0^\infty(M),$ by the integration by parts we have
\begin{align}\label{integbypart}
\begin{split}
&\int_M f^{-\frac{\gamma r}{2}} u^{r-1}\langle \del f,\del u\rangle dV \\
=& -\frac{1}{r}\int_M f^{-\frac{\gamma r}{2}}u^r  \laplace f dV+\frac{\gamma r}{2r}\int_M f^{-\frac{\gamma r}{2} -1} |\del f|^2 u^r dV\\
=& -\frac{1}{r}\int_M f^{-\frac{\gamma r}{2}} u^r \left(\frac{n}{2}-R\right)dV +\frac{\gamma r}{2r}\int_M f^{-\frac{\gamma r}{2}-1}(f-R)u^rdV\\
=&-\frac{n-\gamma r}{2r}\int_M f^{-\frac{\gamma r}{2}}u^r dV+ \frac{1}{r}\int_M f^{-\frac{\gamma r}{2}} u^r R dV -\frac{\gamma}{2}\int_M f^{-\frac{\gamma r}{2}-1} u^r RdV,
\end{split}
\end{align}
where we have used the fact that $\laplace f=\frac{n}{2}-R$ and $f=|\del f|^2+R.$
Then by H\"older's inequality and Cauchy-Schwarz's inequality, there holds
\begin{align*}
&\left(\int_M f^{-\frac{\alpha p}{2}}|\del u|^pdV\right)^\frac{1}{p} \left(\int_M f^{-\frac{(\gamma r-\alpha)p}{2(p-1)}}|\del f|^\frac{p}{p-1}u^\frac{p(r-1)}{p-1}\right)^\frac{p-1}{p}\\
=&\left(\int_M f^{-\frac{\alpha p}{2}}|\del u|^pdV\right)^\frac{1}{p} \left(\int_M f^{-\frac{\beta}{2}}\left(1-\frac{R}{f}\right)^\frac{p}{2(p-1)}u^\frac{p(r-1)}{p-1}\right)^\frac{p-1}{p}\\
\geq & \left| \frac{n-\gamma r}{2r}\int_M f^{-\frac{\gamma r}{2}}u^r dV- \frac{1}{r}\int_M f^{-\frac{\gamma r}{2}} u^r R dV +\frac{\gamma}{2}\int_M f^{-\frac{\gamma r}{2}-1} u^r RdV\right|.
\end{align*}
When $\alpha=0,\gamma r={\widetilde q},r=\widetilde p,p=2,\beta=2\widetilde q-2,$ where $0<\widetilde q<2<\widetilde p,\frac{2(\widetilde p-\widetilde q)}{\widetilde p-2}>n>2,$ the above inequality becomes
\begin{align*}
&\int_M f^{-\tq}|\nabla f|^2 u^{2\tp-2}dV \int_M |\nabla u|^2 dV \\
\geq & \left(\frac{n-\tq}{2\tp}\int_M f^{-\frac{\tq}{2}}u^{\tp} dV - \frac{1}{\tp}\int_M f^{-\frac{\tq}{2}}u^{\tp} R dV+\frac{\tq}{2\tp}\int_M f^{-\frac{\tq}{2}-1}u^{\tp}p R dV
\right)^2, 
\end{align*}
which is equivalent to
  \begin{align*}
\begin{split}
&\int_M |\nabla u|^2 dV \int_M f^{-\tq+1}u^{2\tp-2}dV \\
\geq & \left(\frac{n-\tq}{2\tp}\int_M f^{-\frac{\tq}{2}}u^{\tp} dV - \frac{1}{\tp}\int_M f^{-\frac{\tq}{2}}u^{\tp} R dV+\frac{\tq}{2\tp}\int_M f^{-\frac{\tq}{2}-1}u^{\tp} R dV
\right)^2\\
& + \int_M f^{-\tq}u^{2\tp-2}RdV \int_M |\nabla u|^2 dV.
\end{split}
\end{align*}
 Let $u_\lambda=(\lambda+f^{-\frac{\tq}{2}+1})^\frac{1}{2-\tp}.$ Since $M$ is complete, one can have such $u_\lambda\in \overline{C^\infty_0(M)}.$ Now we put $u_\lambda$ into (\ref{integbypart}), and then have
\begin{align}\label{special}
\begin{split}
&\frac{2-\tq}{2(2-\tp)}\int_M f^{-\tq} |\del f|^2 u_\lambda^{2\tp-2}dV\\ =&-\frac{n-\tq}{2\tp}\int_M f^{-\frac{\tq}{2}} u_\lambda^{\tp} dV + \frac{1}{\tp}\int_M f^{-\frac{\tq}{2}} u_\lambda^{\tp} RdV -\frac{\tq}{2\tp}\int_M f^{-\frac{\tq}{2}-1}u_\lambda^{\tp} RdV.
\end{split}
\end{align}
On the other hand, a direct computation gives 
\begin{align*}
\int_M |\del u_\lambda|^2 dV \int_M f^{-\tq+1}u_\lambda^{2\tp-2}dV= \left(\frac{2-\tq}{2(2-\tp)}\int_M f^{-\tq} |\del f|^2 u_\lambda^{2\tp-2} dV\right)^2.
\end{align*}
Thus the equality of (\ref{special}) holds for such $u_\lambda.$
\endproof

\noindent\textbf{Proof of Corollary 1.1:}
Assume $M$ is non-compact with $R\geq \frac{n}{2}.$
To prove the corollary, one can use some special $u$ according to different parameters. Here, for simplicity, we just use  $u_\lambda=(\lambda+f^{-\frac{\tq}{2}})^\frac{1}{2-\tp}$ that is given in Theorem \ref{pigola}.
By the assumption $R\geq \frac{n}{2}$ and the ranges of $\tp,\tq$, we have that the LHS of (\ref{special}) is non-positive, while the RHS is equal to
\begin{align*}
0\geq LHS= RHS &\geq \frac{\tq}{2\tp}\int_M f^{-\frac{\tq}{2}}u_\lambda^{\tp} dV -\frac{\tq}{2\tp}\int_M f^{-\frac{\tq}{2}-1}u_\lambda^{\tp} RdV\\
& = \frac{\tq}{2\tp}\int_M f^{-\frac{\tq}{2}-1}u_\lambda^{\tp} |\del f|^2dV\\
&>0,
\end{align*}
which gives a contradiction. Thus $f$ has to be a constant, which means that $M$ is Einstein. By Myers's theorem (\cite{Myers}) one has that $M$ has to be compact, which contradicts to our assumption.  
\endproof
\begin{remark}
Note that the choice of $u_\lambda$ in Theorem \ref{pigola} is, in fact, based on the extremal functions of CKN inequalities on $\mathbb R^n$ for different ranges of the parameters (see e.g. \cite{Nguyen}). When the parameters additionally satisfy $r=p$ and $1+\alpha-\frac{\beta}{p}>0,$ we could choose $u_\lambda=e^{-\lambda f^\frac{1+\alpha-\frac{\beta}{p}}{2}},$ which will also give a proof of the final assertion in Theorem \ref{pigola}. In particular, when $\alpha=0,r=p=2,\beta=-2,1+\alpha-\frac{\beta}{p}=1,$ this CKN inequality is reduced to HPW uncertainty principle, and the proof of the final assertion is essentially the same as that in \cite[Theorem 3]{PRS}.
\end{remark}

\noindent\textbf{Proof of Theorem \ref{CKN-thm}}
	Since $M$ is complete, 
	%$$\overline{C_0^\infty(M)}^{W^{1,2}} = \overline{C^\infty(M)\cap W^{1,2}(M)}$$ (see \cite{Hebey}). 
	we can apply (\ref{CKN1}) to $u_\lambda = e^{-\lambda f^\frac{1+\alpha-\frac{\beta}{p}}{2}}$ by a simple approximation procedure (see \cite[Page 49, Theorem 3.1]{Hebey}).
	Consequently, we have
	\begin{align*}
	&\int_M \frac{e^{-p\lambda f^\frac{1+\alpha-\frac{\beta}{p}}{2}}}{f^\frac{\gamma p}{2}}dV\\
	\leq & \frac{2p}{n-\gamma p}\frac{\lambda (1+\alpha-\frac{\beta}{p})}{2}\left(\int_M f^{-\frac{p+\beta}{2}}|\del f|^pe^{-p\lambda f^\frac{1+\alpha-\frac{\beta}{p}}{2}}dV\right)^\frac{1}{p}  \left(\int_M e^{-p\lambda f^\frac{1+\alpha-\frac{\beta}{p}}{2}}f^{-\frac{\beta}{2}}dV\right)^\frac{p-1}{p},
	\end{align*}
	which implies that
	\begin{align*}
	&\int_M \frac{e^{-p\lambda f^\frac{1+\alpha-\frac{\beta}{p}}{2}}}{f^\frac{\gamma p}{2}}dV\\
	\leq &\frac{p\lambda (1+\alpha-\frac{\beta}{p})}{n-\gamma p}\left(\int_M f^{-\frac{\beta}{2}}\left(1-\frac{R}{f}\right)^\frac{p}{2}e^{-p\lambda f^\frac{1+\alpha-\frac{\beta}{p}}{2}}dV\right)^\frac{1}{p}  \left(\int_M e^{-p\lambda f^\frac{1+\alpha-\frac{\beta}{p}}{2}}f^{-\frac{\beta}{2}}dV\right)^\frac{p-1}{p},
	\end{align*}
	Define $F(\lambda) =\int_M  u_\lambda^p f^{-\frac{\gamma p}{2}}dV= \int_M e^{-p\lambda f^{\frac{1+\alpha-\frac{\beta}{p}}{2}-\frac{\gamma p}{2}}} dV.$ Then $F^\prime (\lambda)= -p\int_M f^{\frac{1+\alpha-\frac{\beta}{p}}{2}-\frac{\gamma p}{2}}u_\lambda^p dV$.
	A direct computation shows that 
	\begin{align*}
	F^\prime (\lambda)=& -p\int_M f^{\frac{1+\alpha-\frac{\beta}{p}-\gamma p}{2}}e^{-p\lambda f^{\frac{1+\alpha-\frac{\beta}{p}}{2}}}dV\\
	=&-p\int_M f^{\frac{1+\alpha-\frac{\beta}{p}-\gamma p}{2}-1}\langle \del f,\del f\rangle e^{-p\lambda f^{\frac{1+\alpha-\frac{\beta}{p}}{2}}}dV -p\int_M f^{\frac{1+\alpha-\frac{\beta}{p}-\gamma p}{2}-1}Re^{-p\lambda f^{\frac{1+\alpha-\frac{\beta}{p}}{2}}}dV\\
	=&\frac{2}{\lambda (1+\alpha-\frac{\beta}{p})}\int_M f^{\frac{-\gamma p}{2}}\langle \del f,\del u_\lambda^p\rangle dV -p\int_M f^{\frac{1+\alpha-\frac{\beta}{p}-\gamma p}{2}-1}Ru_\lambda^pdV\\
	=&-\frac{2}{\lambda (1+\alpha-\frac{\beta}{p})}\int_M f^{-\frac{\gamma p}{2}} u_\lambda^p\laplace f dV +\frac{\gamma p}{\lambda (1+\alpha-\frac{\beta}{p})} \int_M|\del f|^2 u^p_\lambda f^{-\frac{\gamma p}{2}-1}dV\\
	&-p\int_M f^{\frac{1+\alpha-\frac{\beta}{p}-\gamma p}{2}-1}Ru_\lambda^pdV\\
	=&-\frac{2}{\lambda (1+\alpha-\frac{\beta}{p})}\int_M f^{-\frac{\gamma p}{2}} u_\lambda^p (\frac{n}{2}-R) dV +\frac{\gamma p}{\lambda (1+\alpha-\frac{\beta}{p})} \int_M (f-R) u^p_\lambda f^{-\frac{\gamma p}{2}-1}dV\\
	&-p\int_M f^{\frac{1+\alpha-\frac{\beta}{p}-\gamma p}{2}-1}Ru_\lambda^pdV\\
	=&-\frac{n-\gamma p}{\lambda(1+\alpha-\frac{\beta}{p})}\int_M f^{-\frac{\gamma p}{2}} u_\lambda^p  dV +\frac{2}{\lambda(1+\alpha-\frac{\beta}{p})}\int_M f^{-\frac{\gamma p}{2}} Ru_\lambda^p  dV\\
	&-\frac{\gamma p}{\lambda (1+\alpha-\frac{\beta}{p})} \int_M R u^p_\lambda f^{-\frac{\gamma p}{2}-1}dV-p\int_M f^{\frac{1+\alpha-\frac{\beta}{p}-\gamma p}{2}-1}Ru_\lambda^pdV,
	\end{align*}
	which means that
	\begin{align*}
	-\frac{F^\prime(\lambda)}{F(\lambda)} =& \frac{n-\gamma p}{\lambda(1+\alpha-\frac{\beta}{p})} - \frac{1}{F(\lambda)}\left( \frac{2}{\lambda(1+\alpha-\frac{\beta}{p})}\int_M f^{-\frac{\gamma p}{2}} Ru_\lambda^p  dV\right.\\
	&\left. -\frac{\gamma p}{\lambda (1+\alpha-\frac{\beta}{p})} \int_M R u^p_\lambda f^{-\frac{\gamma p}{2}-1}dV-p\int_M f^{\frac{1+\alpha-\frac{\beta}{p}-\gamma p}{2}-1}Ru_\lambda^pdV\right)\\
	=&\frac{n-\gamma p}{\lambda(1+\alpha-\frac{\beta}{p})}-\frac{1}{F(\lambda)}\left(I+II+III\right),
	\end{align*}
	where we have used the fact that $\laplace f=\frac{n}{2}-R$ and $f=|\del f|^2+R.$
	Consequently, we have
	\begin{align*}
	-\frac{F^\prime(\lambda)}{F(\lambda)}+\frac{1}{F(\lambda)}(I+II+III)\leq \left(\frac{p\int_M f^{-\frac{\beta}{2}}(1-\frac{R}{f})^\frac{p}{2}u_\lambda^pdV}{F(\lambda)}\right)^\frac{1}{p} \left(-\frac{F^\prime(\lambda)}{F(\lambda)}\right)^{1-\frac{1}{p}},
	\end{align*} 
	and then,
	\begin{align*}
	\frac{1}{F(\lambda)}(I+II+III)\leq \left(\left(\frac{p\int_M f^{-\frac{\beta}{2}}(1-\frac{R}{f})^\frac{p}{2}u_\lambda^pdV}{-F^\prime(\lambda)}\right)^\frac{1}{p}-1\right) \left(-\frac{F^\prime(\lambda)}{F(\lambda)}\right). 
	\end{align*} 
	Obviously,   
	\begin{align*}
	\left(\frac{p\int_M f^{-\frac{\beta}{2}}(1-\frac{R}{f})^\frac{p}{2}u_\lambda^pdV}{-F^\prime(\lambda)}\right)^\frac{1}{p}-1<0
	\end{align*}
	if $R\neq 0.$
	Thus, to get a contradiction, it suffices to show that
	\begin{align*}
	I+II+III\geq 0
	\end{align*} 
	or more exactly,
	\begin{align*}
	I+II+III\geq -F^\prime(\lambda)\left(\left(\frac{p\int_M f^{-\frac{\beta}{2}}(1-\frac{R}{f})^\frac{p}{2}u_\lambda^pdV}{-F^\prime(\lambda)}\right)^\frac{1}{p}-1\right). 
	\end{align*}
	
	Since $\gamma\leq 0,0<1+\alpha-\frac{\beta}{p}\leq 2,$ we have
	\begin{align*}
	III&= -p\int_M f^{\frac{1+\alpha-\frac{\beta}{p}}{2}-1-\frac{\gamma p}{2}} Ru_\lambda^p dV\\
	&\geq -p \epsilon^{\frac{1+\alpha-\frac{\beta}{p}}{4}-\frac{1}{2}}\int_M f^{-\frac{\gamma p}{2}}Ru^p_\lambda dV,
	\end{align*}
	where we have used $f^2>fR>\epsilon>0$ (see \cite{Chow-Lu-Yang}).
	Then we could choose a $\lambda$ such that $I+III>0,$ and consequently, $I+II+III>0.$

\endproof

\end{document}